\input amstex
\documentstyle{amsppt}
\magnification=\magstep1
\baselineskip=18pt
\hsize=6truein
\vsize=8truein

\topmatter
\title On the renormalized volumes for conformally compact
Einstein manifolds
\endtitle
\author Sun-Yung A. Chang\footnote{Princeton University, Dept of
Math, Princeton, NJ 08544, supported by NSF Grant DMS--0245266.},
Jie Qing\footnote{University of California, Dept of Math, Santa
Cruz, CA 95064, supported in part by NSF Grant DMS--0402294.} and
Paul Yang\footnote{Princeton University, Dept of Math, Princeton, NJ
08544, supported by NSF Grant DMS--0245266.}
\endauthor

\leftheadtext{Renormalized volume} \rightheadtext{Chang, Qing and
Yang}

\abstract We study the renormalized volume of a conformally compact
Einstein manifold. In even dimensions, we derive the analogue of the
Chern-Gauss-Bonnet formula incorporating the renormalized volume.
When the dimension is odd, we relate the renormalized volume to the
conformal primitive of the $Q$-curvature. We also show how all the
global information come from the scattering.

\endabstract
\endtopmatter

\document

\head 0. Introduction \endhead

Recently, there is a series of work ([GZ],[FG-2] and [FH]) exploring
the connection between scattering theory on asymptotically
hyperbolic manifolds, the Q-curvature in conformal geometry and the
"renormalized volume" of conformally compact Einstein manifolds. In
particular, in [FG-2], a notion of Q-curvature was introduced for an
odd-dimensional manifold as the boundary of a conformally compact
Einstein manifold of even-dimension. In this note, in section 2
below,  we will clarify the relation between the work of [FG-2] and
the notion of Q-curvature in our earlier work [CQ] for the special
case when the manifold is of dimension three. We then explore this
relation to give a different proof of a result of Anderson [A]
writing the Chern-Gauss-Bonnet formula for conformally compact
Einstein 4-manifold with the renormalized volume. Our proof makes
use of the special exhaustion function introduced in [FG-2] that
yields remarkable simplification in computing the $Q$ curvature. In
section 3, using some recent result of Alexakis on Q-curvature, we
generalize the Chern-Gauss-Bonnet formula involving the renormalized
volume to all even dimensional conformally compact Einstein
manifolds. The formula includes as special case the formula of
Epstein (appendix A in [E]) for conformally compact hyperbolic
manifolds. The formula also shows that the renormalized volume is a
conformal invariant of the conformally compact structure when the
dimension is odd. Finally in section 4, we obtain a formula similar
to that in [FG-2] expressing the renormalized volume of a odd
dimensional conformally compact Einstein manifold as the conformal
primitive of the Q curvature, and in terms of the data of the
scattering matrix.

We would like to mention that there are some interesting recent
works [Ab1] [Ab2] of Albin on the renormalized index theorem.

\head 1. Conformally compact Einstein manifolds and renormalized
volumes \endhead

In this section, we will first recall some basic definitions and facts of
conformally compact Einstein manifolds. We then state the main result in [FG-2].

Given a smooth manifold $X^{n+1}$ of dimension $n+1$ with smooth
boundary $\partial X = M^n$. Let $x$ be a defining function for
$M^n$ in $X^{n+1}$ as follows:
$$
\aligned x>0 \ \ & \text{in $X$};  \\
x=0 \ \ & \text{on $M$} ;  \\
dx \neq 0 \ \ & \text{on $M$}.
\endaligned
$$
A Riemannian metric $g$ on $X$ is conformally
compact if $(X, x^2g)$ is a compact Riemannian manifold with
boundary. Conformally compact manifold $(X^{n+1}, g)$ carries
a well-defined  conformal structure $[\hat g]$ on the boundary $M^n$;
where each $\hat g$ is the restriction of $x^2g$ for a defining
function $x$. We call $(M^n, [\hat g])$ the conformal infinity of
the conformally compact manifold $(X^{n+1}, g)$. If, in addition,
$g$ satisfies $Ric_{g} = -ng$, where
$Ric_g$ denotes the Ricci tensor of
the metric $g$, then we call $(X^{n+1}, g)$ a
conformally compact Einstein manifold.

A conformally compact metric is said to be asymptotically
hyperbolic if its sectional curvature approach $-1$ at $\partial
X=M$. It was shown ([FG-1], [GL]) that if $g$ is an asymptotically hyperbolic metric
on $X$, then a choice of metric $\hat g$ in $[\hat g]$ on $M$
uniquely determines a defining function $x$ near the boundary $M$
and an identification of a neighborhood of $M$ in $X$ with
$M\times (0, \epsilon)$ such that $g$ has the normal form
$$
g = x^{-2} (dx^2 + g_x) \tag 1.1
$$
where $g_x$ is a 1-parameter family of metrics on $M$.

As a conformally compact Einstein metric $g$ is clearly
asymptotically hyperbolic, we have, as computed in [G-1] by Graham,
$$
g_x = \hat g + g^{(2)}x^2 + (\text{even powers of $x$}) +
g^{(n-1)}x^{n-1} + g^{(n)}x^n + \cdots, \tag 1.2
$$
when $n$ is odd, and
$$
g_x = \hat g + g^{(2)}x^2 + (\text{even powers of $x$}) +
g^{(n)}x^n + h x^n\log x + \cdots, \tag 1.3
$$
when $n$ is even. Where $\hat g = x^2 g|_{x=0}$, $g^{(2i)}$ are
determined by $\hat g$ for $2i < n$. The trace part of $g^{(n)}$
is zero when $n$ is odd; the trace part of $g^{(n)}$ is determined
by $\hat g$ and $h$ is traceless and determined
by $\hat g$ too when n is even.

As a realization of the holography principle proposed in physics,
one considers the asymptotic of the volume of a conformally
compact Einstein manifold $(X^{n+1}, g)$. Namely, if denote by $x$
the defining function associated with a choice of a metric $\hat g
\in [\hat g]$, we have
$$
\text{Vol}_g(\{x > \epsilon\}) = c_0 \epsilon^{-n} +
c_2\epsilon^{-n+2} + \cdots + c_{n-1} \epsilon^{-1} + V + o(1)
\tag 1.4
$$
for $n$ odd, and
$$
\text{Vol}_g(\{x > \epsilon\}) = c_0 \epsilon^{-n} +
c_2\epsilon^{-n+2} + \cdots + c_{n-2}\epsilon^{-2} + L \log \frac
1\epsilon + V + o(1) \tag 1.5
$$
for $n$ even. We call the constant term $V$ in all dimensions the
renormalized volume for $(X^{n+1}, g)$. We recall that $V$ in odd dimension and $L$
in even dimension are
independent of the choice $\hat g$ in the class $[\hat g]$ (cf.
[HS] [G-1]).

Based on the work of [GZ], Fefferman and Graham [FG-2] introduced the following
formula to calculate the
renormalized volume $V$ for a conformally compact Einstein
manifold. Here we will quote a special case of their result.
For odd $n$, upon a choice of a special defining
function $x$, one sets
$$ v = - \frac{d}{ds} |_{s=n} \wp (s) 1,
$$
where $\wp (s)$ denotes the Possion operator (see [GZ] or section 4
below for the definition of the operator) on $X^{ n+1}$.
The $v$ solves
$$
-\Delta v = n \quad\text{in $X^{n+1}$}, \tag 1.6
$$
and has the asymptotic behavior
$$
v = \log x + A + B x^n  \tag 1.7
$$
in a neighborhood of $M^n$, where $A, B$ are functions even in $x$, and $A|_{x=0}=0$. Then Fefferman and
Graham in [FG-2] defined
$$
(Q_n)_{(g, \hat g)} = k_n B|_{x=0} \tag 1.8
$$
where $k_n = 2^n\frac {\Gamma(\frac n2)}{\Gamma(-\frac n2)}$.

\proclaim{Theorem 1.1} ([FG-2]) When $n$ is odd,
$$
V (X^{n+1}, g) = \frac 1{k_n} \int_M (Q_n)_{(g,
\hat g)} dv_{\hat g}.
\tag 1.9
$$
\endproclaim

In section 2 and 3 below, we will apply Theorem 1.1 to identify
the renormalized volume as part of the integral in the
Gauss-Bonnet formula for $(X^{n+1}, M^n, g) $ when n is odd. In
section 4, we recall the work in [GZ] and [FG-2] relating the
Q-curvature to data of scattering matrix on the asymptotically
hyperbolic manifold $X^{n+1}$ and derive a similar formula for the
renormalized volume when n is even.

\head 2. Chern-Gauss-Bonnet formula for $n=3$ \endhead

To motivate our discussions in this section we first recall some
works in [CQ]. First we recall that the Paneitz
operator defined on 4-manifold as:
$$
P_{4/2} = \Delta^2 + \delta\{\frac 23 Rg - 2\text{Ric} \}d, \tag
2.1
$$
where $R$ is the scalar curvature, $\text{Ric}$ is the Ricci
curvature. There are two important properties of the Paneitz operator:

$$
(P_{4/2})_{g_{w}} = e^{-4 w} (P_{4/2})_{g}.
$$
$$
(P_{4/2})_{g} w  + (Q_4)_g = (Q_4)_{g_w} e^{4 w}, \tag 2.2
$$
for any smooth function $w$ defined on the 4-manifold, and
where $g_{w} = e^{ 2 w} g$ and where $Q$ is the curvature
function
$$
Q_4 = \frac {1}{6} (- \Delta R + |R|^2 - 3 |Ric|^2).
$$

In [CQ], on a compact Riemannian 4-manifold $(X^4, g)$ with boundary,
a third order boundary operator $P_b$ and a third order
boundary curvature  $T$  were introduced as follows:
$$
(P_b)_{g} = -\frac 12 \frac {\partial}{\partial n} \Delta_g +
\tilde\Delta \frac {\partial}{\partial n} + \frac 23 H
\tilde\Delta + L_{\alpha\beta}\tilde\nabla_\alpha
\tilde\nabla_\beta + \frac 13 \tilde\nabla_\alpha H \cdot
\tilde\nabla_\alpha - (F-\frac 13 R)\frac {\partial}{\partial n}
\tag 2.3
$$
and
$$
T_{g} = \frac 1{12} \frac {\partial R}{\partial n} + \frac 16 RH -
R_{\alpha N \beta N}L_{\alpha\beta} + \frac 19 H^3 - \frac 13
\text{Tr} L^3 - \frac 13 \tilde \Delta H, \tag 2.4
$$
where $\frac {\partial}{\partial n}$ is the outer normal
derivative, $\tilde\Delta$ is the trace of the Hessian of the
metric on the boundary, $\tilde\nabla$ is the derivative in the
boundary, $L$ is the second fundamental form of boundary, $H =
\text{Tr}L$, $F = R_{\alpha N\alpha N}$, and $R$ is the scalar
curvature all with respect to the
metric $g$. $P_b$ and $T$ transform under conformal change of metric
on the boundary of $X^4$ similar to that of $P_{4/2}$ and $Q_4$ on
$X^4$ as follows:
$$
\aligned
(P_b)_{g_{w}} & = e^{-3 w} (P_b)_{g} \\
(P_b)_{g} w + T_g & = T_{g_w} e^{3 w}
\endaligned
\tag 2.5
$$
We remark that $(P_b)_{g}$ and $T_g$ thus defined depend on the metric
$g$ on $(X^4, g )$, and are not intrinsic quantities on the
boundary of  $X^4$.

In [CQ], we have also re-organized the terms in the integrand of Gauss-Bonnet formula for
4-manifolds with boundary into the following form:
$$
\frac 1{8\pi^2}\int_{X^4} (|\Cal W|^2 + Q)dv + \frac 1{4\pi^2}
\int_{\partial X} (\Cal  L + T)d\sigma = \chi(X), \tag 2.6
$$
where
$$
\Cal L = \frac 13 RH - FH + R_{\alpha N \beta N}L_{\alpha\beta} -
R_{\alpha\gamma\beta\gamma}L_{\alpha\beta} +\frac 29 H^3 - H|L|^2
+ \text{Tr}L^3. \tag 2.7
$$
We remark that the Weyl curvature $\Cal W$ is a point-wise conformal invariant term
on the 4-manifold, while $\Cal L$ is a point-wise conformal invariant
term on the boundary of the manifold.

We also remark that when the boundary is totally geodesic, the
expressions of $(P_b)_{g}$ and $T_g$ in (2.3) and (2.4) above become
very simple;
$$
(P_b)_{g} = - \frac 12 \frac {\partial}{\partial n} \Delta_g
\, + \, \tilde\Delta \frac {\partial}{\partial n} -
(F-\frac 13 R)\frac {\partial}{\partial n}, \,\,\,\, \,
T_g = \frac {1}{12} \frac {\partial R } {\partial n}, $$
and in this case $\Cal L$ vanishes.

Given a conformally compact Einstein manifold $(X^{n+1}, g)$ and a
choice of metric $\hat g$ in the conformal infinity $(M, [\hat
g])$, we consider the compactification $(X^{n+1}, e^{2v}g)$, where
$v$ is the function which satisfies (1.6) and (1.7). We observe that

\proclaim{Lemma 2.1} When $n$ is odd, $
(Q_{n+1})_{e^{2v} g} = 0$.
\endproclaim
\demo{Proof} The proof follows an observation made by Graham
([G-2], see also [Br]) that the Paneitz operator $P_{\frac {n+1}2}$
on an Einstein manifold is a polynomial of the Laplacian $\Cal P
(\Delta)$ and the polynomial $\Cal P$ on the Einstein manifold is
the same as the one on the constant curvature space with the constant
the same as the constant of the scalar curvature of the Einstein manifold.
Meanwhile the $Q$-curvature $Q_{n+1}$ of an Einstein manifold is
the same as the one on the constant curvature space. Therefore
$(P_{\frac {n+1}2})_{g} = \Cal P (\Delta_g)$
 if $(P_{\frac{n+1}2})_{g_H}= \Cal P
(\Delta_{g_H})$, and $(Q_{n+1})_{g} =
(Q_{n+1})_{g_H}$, where $(H^{n+1}, g_H)$ is the hyperbolic space.
$$
(P_{\frac {n+1}2})_{g_{H^{n+1}}} = \prod_{l=1}^{\frac {n+1}2}
(-\Delta_{H^{n+1}} - C_l) \tag 2.8
$$
where $C_l = (\frac {n+1}2 +l-1)(\frac {n+1}2 - l)$. Therefore
$$
(P_{\frac{n+1}2})_{g} = \sum_{l=2}^{\frac {n+1}2} (-1)^{\frac {n+1}2
-l} B_l (\Delta_g)^l - (-1)^{\frac
{n-1}2} (n-1)! \Delta_g. \tag 2.9
$$
Meanwhile $(Q_{n+1})_{H^{n+1}} = (-1)^{\frac {n+1}2} n!$. Thus
$$
(Q_{n+1})_g = (-1)^{\frac {n+1}2} n!. \tag 2.10
$$
Thus if $v$ satisfies the equation (1.6), we have
$$
(P_{\frac {n+1}2})_{g} v + (Q_{n+1})_{g} = 0. \tag 2.11
$$
It thus follows from equation (2.2) that $
(Q_{n+1})_{ e^{2v} g} = 0$.
\enddemo

We will now combine the above observation to Theorem 1.1
of [FG-2] to give an alternative proof of a result of Anderson [A]
(Theorem 2.3 below) for conformal compact Einstein 4-manifold $(X^4,g)$.
We first
relate our curvature $T$ to that of $Q_3$ as defined in (1.8).

\proclaim{Lemma 2.2}
$$
T_{e^{2v}g} = 3B|_{x=0} = (Q_3)_{(g, \hat g)}. \tag 2.12
$$
\endproclaim

\demo{Proof} By the scalar curvature equation we have
$$
\frac 1{12} R_{e^{2v}g} = \frac 12 (-\Delta_g
e^v + \frac 16 R_g
e^v)e^{-3v}.
$$
Therefore for $v$ satisfies equation (1.6), we have
$$
\frac 1{12} R_{e^{2v}g} = \frac 12 ((e^{-v})^2 - |\nabla
e^{-v}|^2).
$$
We now apply the asymptotic  expansion of $v$ in (1.7) and write
$$
\aligned e^{-2v} & = \frac 1{x^2} - 2A_2 - 2B_0x + O(x^2) \\
|\nabla e^{-v}|^2 & = \frac 1{x^2} + 2A_2 + 4B_0 x + O(x^2),
\endaligned
$$
where  $A_2$ is the coefficient of $x^2$ of $A$ and $B_0 = B|_{ x=0} $.
We get
$$
T_{e^{2v}g} = - \frac 1{12}\frac {\partial}{\partial x}
R_{e^{2v}g}|_{x=0} =3B_0 = Q_{(g, \hat g)}.
$$
This finishes the proof of the lemma.
\enddemo

\proclaim{Theorem 2.3} [A]  Suppose that $(X^4, g)$ is a conformally
compact Einstein manifold. Then
$$
\frac 1{8\pi^2} \int_{X^4} |\Cal W|_g^2 dv_g + \frac 3{4\pi^2}V (X^4,
g) = \chi(X^4). \tag 2.13
$$
\endproclaim
\demo{Proof} Apply Lemma 2.1 to (2.6), we have
$$
\frac 1{8\pi^2}\int_{X^4} |\Cal W|_{e^{2v} g}
\, dv_{e^{2v} g} + \frac
1{4\pi^2} \int_M (\Cal L + T)_{ (e^{2v}g, \hat
g)}\, dv_{\hat g} = \chi(X^4).
$$

We now observe that as the boundary of $M$ of $X^4$ is umbilical,
$\Cal L_{ (e^{2v}g, \hat g )} = 0$. Apply Lemma 2.2 and Theorem 1.1.
we obtain (2.13) for the metric $e^{2v} g$. We then observe that once
the formula (2.13) holds for the metric $e^{2v} g$, it holds for any
metric $\tilde g \in [g]$ with $(X^{n+1}, \tilde g)$ a conformally
compact manifold as the term of the renormalized volume $V$ is
conformally invariant.
\enddemo

\head 3. Chern-Gauss-Bonnet formula in higher dimensions
\endhead

In higher dimensions when $n=2k +1 >3$ we wish to determine the
analogous formula for the Euler characteristic. We continue to
consider the metric $(X^{n+1},e^{2v}g)$ where $v$ satisfies the
equations (1.6) and (1.7). We will find that the parity conditions
imposed in (1.7) makes it possible to determine the local boundary
invariants of order $n$ for the compact manifold $(X^{n+1},
e^{2v}g)$. According to (1.1) and (1.7) we have the expansion of
the metric $e^{2v} g$.
$$
\aligned e^{2v} g & = H^2 dx^2 + \hat g + c^{(2)}x^2 +  \text{even
powers in $x$} \\
& \quad\quad  + c^{(n-1)}x^{n-1} + (2B_0\hat g + g^{(n)})x^n +
\cdots \endaligned \tag 3.1
$$
where
$$
H = e^{A+Bx^n} = 1 + e_2x^2 +  \text{even powers in $x$} +
e_{n-1}x^{n-1} + B_0x^n + \cdots
$$
and $c^{(2i)}$ for $ 1\leq i \leq (n-1)/2$ are local invariants of
$\hat g$. We remark that it is easy to see that the boundary of
$(X^{n+1}, e^{2v}g)$ is totally geodesic.

\proclaim{Lemma 3.1}
$$
(\partial_x \Delta^{\frac {n-3}2} R)_{e^{2v}g}|_{x=0} = -2nn!B_0.
\tag 3.2
$$
\endproclaim
\demo{Proof} We have
$$
\aligned \Delta_{ e^{2v} g} & = \frac 1{H\sqrt {\det g_x}}\partial_\alpha(H
\sqrt {\det g_x}g^{\alpha\beta}_x \partial_\beta)\\
& = Q^{(2)}_2\partial^2_x + Q^{(1)}_2\partial_x + Q^{(0)}_2
\endaligned
$$
where the coefficients $Q^{(i)}$ have the following properties:
$Q^{(2)}_2$ is a zeroth order differential operator, having an
asymptotic expansion in powers of $x$ in which the first
nonzero odd power term is $x^n$.
$Q^{(1)}_2$ is a zeroth order differential operator, having
an expansion in which
the first nonzero even degree term is $x^{n-1}$.
$Q^{(0)}_2$ is differential operator of order 2 of purely
tangential differentiations with coefficients which have expansion in
powers of $x$ in which the  first nonzero odd term is $x^n$.
Inductively, we see
that, for $k \leq \frac {n-3}2$,
$$
\Delta^k = Q^{(2k)}_{2k}\partial_x^{2k} +
Q^{(2k-1)}_{2k}\partial_x^{2k-1} + \cdots + Q^{(1)}_{2k}\partial_x
+ Q^{(0)}_{2k} \tag 3.3
$$
where $Q^{(i)}_{2k}$ ($i\neq 0$) is a differential operator of
order $2k-i$ of purely tangential differentiations with
coefficients having expansions in powers of $x$ in which the first nonzero even
terms are $x^{n-(2k-i)}$ if $i$ is odd, and the first nonzero odd
terms are $x^{n-(2k-i)}$ if $i$ is even, and $Q^{(0)}_{2k}$ is a
differential operator of order $2k$ of purely tangential
differentiations with coefficients whose expansions in $x$ have
the first nonzero odd terms $x^{n-2k+2}$. Thus
$$
\partial_x\Delta^k = F^{(2k+1)}\partial_x^{2k+1} +
F^{(2k)}\partial_x^{2k} + \cdots + F^{(1)}\partial_x + F^{(0)}
\tag 3.4
$$
where $F^{(2k+1)} = Q^{(2k)}$, $F^{(i)}$ ($ 0<i<2k+1$) is a
differential operator of order $2k-i+1$ of purely tangential
differentiations with coefficients whose expansions in $x$ have
the first nonzero even terms are $x^{n-(2k-i)-1}$ if $i$ is even,
and the first nonzero odd terms are $x^{n-(2k-i)-1}$ if $i$ is
odd, and $F^{(0)}$ is a differential operator of order $2k$ of
purely tangential differentiations with coefficients whose
expansions in $x$ have the first nonzero even terms $x^{n-2k+1}$.

On the other hand, we have
$$
R_{e^{2v}g} = -2 n^2(n-1) B_0 x^{n-2} + \text{even powers of $x$
terms} + o(x^{n-2}). \tag 3.5
$$
Keeping track of the parity, we obtain (3.2).
\enddemo

Next we deal with all other boundary terms, these are
contractions of one or more factors consisting of curvatures, covariant
derivatives of curvatures, except $\partial_x^{n-2}R$ which
is accounted in the above term $\partial_x\Delta^{\frac {n-3}2}R$.
Since $n$ is odd, and $\partial x$ is the normal direction,
each such term must contain at least one $x$
index. In fact, the total number of $x$ indices appearing in each of
such terms must be odd. Thus one finds that each of such terms
always contains a factor which is a covariant derivatives of
curvature and in which $x$ index appears odd number of times. Such
factors, if we insist on taking $\nabla_x$ first, must appear as one
of the following three different types
$$
\nabla_\spadesuit\cdots\nabla_\spadesuit\nabla_x^{2k+1}
R_{\spadesuit\spadesuit\spadesuit\spadesuit} \tag I
$$
where $\spadesuit$ stands for indices other than $x$, in other
words, tangential.
$$
\nabla_\spadesuit\cdots\nabla_\spadesuit\nabla_x^{2k}
R_{x\spadesuit\spadesuit\spadesuit} \tag II
$$
and
$$
\nabla_\spadesuit\cdots\nabla_\spadesuit\nabla_x^{2k-1}
R_{x\spadesuit x\spadesuit}. \tag III
$$
Note that in all three types $1 \leq 2k+1 \leq n-2$. Since the
boundary is totally geodesic, we only need

\proclaim{Lemma 3.2} All three types of boundary terms
$$
\nabla_x^{2k+1} R_{\spadesuit\spadesuit\spadesuit\spadesuit},
\quad \nabla_x^{2k} R_{x\spadesuit\spadesuit\spadesuit}, \quad
\nabla_x^{2k-1} R_{x\spadesuit x\spadesuit}
$$
vanish at the boundary for $1 \leq 2k+1 \leq n-2$.
\endproclaim
\demo{Proof} We consider a point at the boundary and choose a
normal coordinate on the boundary $M^n$ in the special coordinates
for $X^{n+1}$. Recall
$$
\aligned R_{\alpha\beta\gamma\delta} & = \frac 12
(-\partial_\beta\partial_\delta g_{\alpha\gamma} -
\partial_\alpha\partial_\gamma g_{\beta\delta} +
\partial_\beta\partial_\gamma g_{\alpha\delta} +\partial_\alpha\partial_\delta g_{\beta\gamma}) \\
& - g^{\eta\lambda}([\alpha\gamma,\eta][\beta\delta,\lambda] -
[\beta\gamma,\eta][\alpha\delta,\lambda]),
\endaligned
$$
and
$$
\nabla_x T_{\alpha\beta\cdots\delta} = \partial_x
T_{\alpha\beta\cdots\delta} - \Gamma^\lambda_{\alpha \ x}
T_{\lambda\beta\cdots\delta}- \Gamma^\lambda_{\beta \ x}
T_{\alpha\lambda\cdots\delta} - \cdots  - \Gamma^\lambda_{\delta \
x} T_{\alpha\beta\cdots\lambda}
$$
where
$$
\Gamma^\alpha_{\beta\gamma} = g^{\alpha\delta}[\beta\gamma,
\delta]
$$
and
$$
[\alpha\beta,\gamma] = \frac 12 (\partial_\alpha g_{\beta\gamma} +
\partial_\beta g_{\alpha\gamma} - \partial_\gamma g_{\alpha\beta}).
$$
For the simplicity of the notation we will use $g$ to stand for
$e^{2v}g$ if no confusion can arise. Each of the three types is a sum
of products of factors that are of the form:
$$
\partial_\alpha\partial_\beta\cdots\partial_\gamma g_{\lambda\mu}
$$
or
$$
\partial_\alpha\partial_\beta\cdots\partial_\gamma g^{\lambda\mu}.
$$
We claim that each summand must has a factor that is one of the
following
$$
\partial_\spadesuit\cdots\partial_\spadesuit\partial_x^{2k+1}g_{\spadesuit\spadesuit},
$$
$$
\partial_\spadesuit\cdots\partial_\spadesuit\partial_x^{2k-1}g_{xx},
$$
$$
\partial_\spadesuit\cdots\partial_\spadesuit\partial_x^{2k+1}g^{\spadesuit\spadesuit}.
$$
and
$$
\partial_\spadesuit\cdots\partial_\spadesuit\partial_x^{2k-1}g^{xx}.
$$
where  $1\leq 2k+1 \leq n-2$.
To verify the claim, one needs to observe that, in writing
the three types in local coordinates, the number of times the index $x$
appears in each
summand increases only when one sees
$$
\Gamma^x_{\spadesuit \ x} T_{\alpha\beta\cdots\ x \cdots \delta},
$$
where the number of $x$ increases by $2$. Thus, in the end, the
total number of index $x$ in each summand is still odd. Therefore one of
the factors must have odd number of $x$. Finally one observes that
for any individual factor arising here the number of $x$ can not
exceed $n-1$. So the proof is complete.
\enddemo

We now apply above results to derive a formula analogous to that
of Theorem 2.3 for the renormalized volume on $(X^{n+1}, g)$ for
$n=5$. In this case, we recall the formula of Graham [G-1]:
$$
Q_6 = 64 \pi^3 e - \frac 16 J + \frac 1{10} \Delta^2 R +
\text{Div}(T) \tag 3.6
$$
where $e$ is the Euler class density whose integral over a compact
6-manifold gives its Euler number,
$$
J  = - 3 I + 7W_{ij}^{\quad ab}W_{ab}^{\quad pq}W_{pq}^{\quad ij}
+ 4 W_{ijkl}W^{iakb}W^{j \ l}_{\ a \ b},
$$
$$
I = |V|^2 - 16 W^i_{\ jkl}C_{jkl,i} + 16 P^{im}W_{ijkl}W_m^{\ jkl}
+ 16 |C|^2,
$$
$$
V_{ijklm} = W_{ijkl,m} + g_{im}C_{jkl} - g_{jm}C_{ikl} +
g_{km}C_{lij} - g_{lm}C_{kij},
$$
$$
C_{ijk} = P_{ij,k} - P_{ik,j},
$$
$$
P_{ij} = \frac 14(R_{ij} - \frac R{10}g_{ij}),
$$
$W_{ijkl}$ is the Weyl curvature, $R_{ij}$ is the Ricci curvature,
$R$ is the scalar curvature, and $\text{Div}(T)$ are divergence
terms other than $\Delta^2 R$. $Q_6$ in this form is organized so
that it is a sum of three types of terms: Euler class density,
local conformal invariants, and divergence terms.
G
\proclaim{Theorem 3.3} For $n=5$, we have
$$
\chi (X^6) =  \frac 1{128\pi^3} \int_{X^6}
(\Cal J)_g  dv_g   - \frac
{15}{8 \pi^3}V (X^6, g),  \tag 3.7
$$
where
$$
\Cal J = - |\nabla  W|^2 + 8 |W|^2 + \frac 73 W_{ij}^{\quad ab}
W_{ab}^{\quad pq}W_{pq}^{\quad ij} + \frac 43 W_{ijkl}
W^{iakb}W^{j \ l}_{\ a \ b}
$$
\endproclaim
\demo{Proof} As a consequence of above lemmas, we have
$$
\chi (X^6) =  \frac 1{128\pi^3} \int_{X^6} (\frac 13
 J_{ e^{2v} g} ) dv_{e^{2v}
g} - \frac {15}{8 \pi^3}V (X^6, g). \tag 3.8
$$
Since $J$ is a local conformal
invariant and $g$ is an Einstein metric, we obtain (3.7) directly from (3.8).
\enddemo

We can find a general formula for all higher dimensions. We recall
a recent result of S. Alexakis [Al].

\proclaim{Theorem 3.4} [Alexakis] On any compact
Einstein m-dimensional
manifold with m even, we have
$$
Q_{m} = a_{m} e + \Cal J +
\text{Div}(T_m). \tag 3.9
$$
where $e$ is the Euler class density, $\Cal J$ is a conformal
invariant, and $\text{Div}(T_m)$ is a divergence term and $a_m$ is some dimensional constant.
\endproclaim

We also recall a fact we learned from Tom Branson [Br]:

\proclaim {Proposition 3.5} On any compact m-dimensional manifold
for m even, suppose $Q_m$ is {\it the} curvature in the
construction of [GJMS] and [GZ], [FG-2], then
$$
Q_m = b_m {\Delta}^{ \frac {m-2}{2} } R + \text{lower order
terms}, \tag 3.10
$$
where
$$
b_m = (-1)^{ \frac {m-2}2} \frac{ 2^{m-1} (\frac {m} 2)!\Gamma
(\frac {m+1}2)} {\sqrt \pi (m-1) m!}.
$$
\endproclaim


\proclaim{Theorem  3.6} When $n$ is odd, we have
$$
\int_{X^{n+1}} (\Cal W_{n+1})_g dv_g + (-1)^{\frac
{n+1}2} \frac {\Gamma(\frac {n+2}2)}{\pi^{\frac {n+2}2}} V
(X^{n+1}, g) = \chi (X^{n+1}) \tag 3.11
$$
for some curvature invariant $\Cal W_{n+1}$, which is a sum of
contractions of Weyl curvatures and/or its covariant derivatives
in an Einstein metric.
\endproclaim

\demo{Proof} We first establish that equation (3.9) remains valid on a conformally Einstein manifold $(X^{n+1}, g)$.
Let $g_w=e^{2w}g$ be such a
metric, then it follows from the Paneitz equation that for
$m = n+1$,
$$
\aligned (Q_m)_{g_w} e^{v} &= (P_m)_g v \, + \,  (Q_m)_ g \\
&=a_m e_g + {\Cal J}_g + Div(T')\\
&=a_m e_{g_w} + {\Cal J}_{g_w} + Div(T'')
\endaligned
\tag 3.9'
$$
the second equation follows from the fact that the Paneitz operator
$P_m$ is a
divergence and Theorem 3.4. The third equation follows from the fact that the
Pfaffians of any two Riemannian metrics on the same manifold
differs by a divergence term and ${\Cal J}$ is a conforaml invariant.

In order to apply this formula, we need to observe that the
leading order term $\Delta ^{\frac{m-2}{2}} R$ in formula (3.10)
cannot appear in the conformally invariant term $\Cal J$. In order to
see this, we first recall that the $\Cal J$ is a linear
combination of terms of the form $Tr( \nabla^{I_1}{\Cal
R}\bigotimes \nabla^{I_2}{\Cal R}...\bigotimes \nabla^{I_k}{\Cal
R} ) $ of weight $m$ where $Tr$ denotes a suitably chosen pairwise
contraction over all the indices. Observe that the conformal
variation $\delta_w ( \Delta ^{\frac{m-2}{2}})  R$, where $\delta_w$ denotes the variation of the metric $g$ to $g_w$
 is of the form $
\Delta ^{\frac{m}{2}} w  + \text{lower order terms}$. Thus if
$\Delta ^{\frac{m-2}{2}} R$ does appear as a term in $\Cal J$, its
conformal variation must be cancelled by the conformal variations
of the other terms in the linear combination, but it is clear that
the conformal variations of the other possibilities
of the curvature $\Cal R$  other than the scalar curvature $R$ cannot have
order m in the number of derivatives of $w$ and of the form
$\Delta ^{\frac{m}{2}} w$.

We can now apply the formula (3.9') to the metric $g_v = e^{2v}g$ where $v$ is as in Lemma 2.1., thus by Lemma 2.1 the left hand side of (3.9') is
identically zero, and we find
$$
a_m \chi(X^{n+1})= \int_{X^{n+1}} ({\Cal J}_{g_v}  - Div(T''))
 dv_{g_v}.
$$

Among the divergence terms in $Div(T'')$, only the leading
order term $ b_m \Delta^{ \frac {m-2}{2}} R$ has a non-zero contribution according to Lemma 3.2. The
computation in Lemma 3.1 determines the precise contribution of
this term as a mutiple of the renormalized volume. We also note that as $g$ is an Einstein metric, we may assume that
the terms which appear in the conformal invariant
${\Cal J} $ are contractions of the Weyl curvature together with its
covariant derivatives. We have thus finished the proof of Theorem 3.6.
\enddemo

 \proclaim{Corollary 3.5} When $(X^{n+1}, g)$ is conformally
compact hyperbolic, we have
$$
V(X^{n+1}, g) = \frac {(-1)^{\frac {n+1}2} \pi^{\frac
{n+2}2}}{\Gamma(\frac {n+2}2)} \chi(X). \tag 3.12
$$
\endproclaim

One may compare (3.12) to a formula for renormalized volume given
by Epstein in [E], where he has
$$
V(X^{n+1}, g) = \frac {(-1)^m 2^{2m} m!}{(2m)!} \chi(X) \tag 3.13
$$
for $n= 2m-1$.

\head 4. Scattering theory and the renormalized volume
\endhead

We now recall the connection between the renormalized volume and
scattering theory introduced in [GZ]. Suppose that $(X^{n+1}, g)$
is a conformally compact Einstein manifold and $(M^n, [\hat g])$
is its associated conformal infinity. And suppose that $x$ is a
defining function associated with a choice of metric $\hat g\in
[\hat g]$ on $M$ as before. One considers the asymptotic Dirichlet
problem at infinity for the Poisson equation
$$
(-\Delta_g - s(n-s)) u = 0. \tag 4.1
$$
Based on earlier works on the resolvents, Graham and Zworski in
[GZ] proved that there is a meromorphic family of solutions $u (s)
= \wp(s)f$ such that
$$
\wp(s) f = F x^{n-s}  + G x^s  \ \  \text{if $s \notin n/2 + N$}
\tag 4.2
$$
where $F, G, \in C^\infty(X)$, $F|_M = f$, and $F, G$ mod
$O(x^n)$ are even in $x$. Also if
$n/2 - s$ is not an integer, then $G|_M$ is globally determined by
$f$ and $g$. The scattering operator is defined as:
$$
S(s)f = G|_M. \tag 4.3
$$


Thus the function $v$ satisfying (1.6) and (1.7) studied in [FG-2] is defined as:
$$
v = - \frac {d}{ds}|_{s=n} \wp(s)1. \tag 4.4
$$
Therefore when $n$ is odd, we may rewrite (1.8) in Theorem 1.1 as
$$
V (X^{n+1}, g) = - \int_M \, (\frac{d}{d
s}|_{s=n} S(s)1) \,  dv_{\hat g}.
\tag 4.5
$$

We will now point out that formula similar to that (4.5) holds also
when $n$ is even. To do so, we first establish some notations.

When $n$ is even and $(X^{n+1}, g) $ conformal Einstein. For each  $ s<n$ and close to $n$, we consider the
solution of the Possion equation $ \wp(s) 1 = u(s)$ as in (4.1).
Then
$$
u(s) = x^{n-s} F(x,s) + x^s G(x,s) \tag 4.6
$$
for functions $F, G$ which are even in x. We denote the asymptotic
expansion of $ F$ near boundary as
$$
F(x,s) = 1 + a_2(s) x^2 + a_4 (s) x^4 + \cdot \,\, +a_n(s) x^n + \cdot.
\tag 4.7
$$
Denote $ a_k' = \frac {d}{ds}|_{ s=n}  a_k (s) $. We have the following
formula.

\proclaim{Theorem 4.1} Suppose that $(X^{n+1}, g)$ is a
conformally compact Einstein manifold with even $n$. For a choice
of metric $\hat g\in [\hat g]$ of the conformal infinity $(M,
[\hat g])$, the renormalized volume is
$$
\aligned & V(X^{n+1}, g, \hat g)  = - \int_M
(\frac {d}{d
s}|_{s=n} S(s)1 ) \,  dv_{\hat g} \\ & \quad - \frac 1n \int_M 2a_2'
v^{(n-2)}dv_{\hat g} - \cdots - \frac 1n \int_M
(n-2)a_{n-2}'v^{(2)}dv_{\hat g} - \int_M a_n'dv_{\hat g} .
\endaligned
\tag 4.8
$$
\endproclaim
\demo{Proof}
Using the notations as set in (4.6) and (4.7), we have
$$
v = - \frac{d}{d s}|_{s=n} u(s) = F(x,n) \log x  - F' - x^n G' + G (x,
n) x^n\log x.
$$
Since
$$
F(x,n) = 1 - c_{\frac n2} Q_n x^n + O(x ^{ n+1}),
$$
and
$$ G(x, n)|M = c_{\frac n2} Q_n =  \lim_{s \to n} S(s)1,
$$
where
 $c_k = (-1)^k (2^{2k}k!(k-1)!)^{-1}$,
and
$$
F' = \frac {d F}{d s}|_{s=n}, \quad \text{and} \ G' = \frac {d
G}{d s}|_{s=n}.
$$
We get
$$
v \, = \,  \log x -F' - x^n G' - 2 c_{\frac n2} Q_n x^n \log x + O( x^{n+1} \log x).
$$
Recall the expansion of the volume element
$$
\aligned
dv_{g_\epsilon} & = \sqrt {\frac {\det g_\epsilon}{\det
\hat g}} dv_{\hat g}  = \sqrt {\det (\hat g^{-1}g_{\epsilon})} dv_{\hat g}\\
& = (1 + v^{(2)} \epsilon^2 + v^{(4)} \epsilon^4 + \cdots)dv_{\hat
g}
\endaligned
$$
where $v^{(2)}, \cdots, v^{(n)}$ are determined by $\hat g$. We
have
$$
\text{vol}(\{x >\epsilon\}) = \frac 1n \int_{x > \epsilon} -\Delta
v dv_g = \frac 1n \int_{x = \epsilon} -\frac {\partial v}{\partial
n}d\sigma_{g_\epsilon} = \frac 1n \epsilon^{-n+1} \int_M \frac
{\partial v}{\partial x}|_{x = \epsilon}dv_{g_\epsilon}
$$
where
$$
\aligned
\frac {\partial v}{\partial x} \, = \, & \frac 1x - \cdots - na_{n}'
x^{n-1} - n (\frac {d }{d s}|_{s=n} S(s)1) x^{n-1} \\ & - 2c_{\frac
n2}Q_nx^{n-1} - 2n c_{\frac n2} Q_n x^{n-1} \log x +
o(x^{n-1}). \\
\endaligned
$$
Thus
$$
\aligned  \text{vol}(\{x > \epsilon\})
& = \cdots + \frac 1n \int_M (v^{(n)} - 2a_2' v^{(n-2)} - \cdots -
(n-2)v^{(2)}
a_{n-2}' \\
& \quad\quad - na_n' - n \frac {d }{d s}|_{s=n} S(s)1 - 2c_{\frac
n2}Q_n )dv_{\hat g} + \cdots
\endaligned
$$
and the renormalized volume for $(X^{n+1}, g)$ is
$$
\aligned V(X^{n+1}, g, \hat g) & = - \frac 1n \int_M 2a_2'
v^{(n-2)}dv_{\hat g} - \cdots -
\frac 1n \int_M (n-2)a_{n-2}'v^{(2)}dv_{\hat g} \\
& \quad\quad - \int_M a_n'dv_{\hat g} - \int_M
\, (\frac {d }{d
s} |_{s=n} S(s)1) \,  dv_{\hat g}.
\endaligned
$$
We have thus established the formula (4.8).
\enddemo

We remark that the coefficients $a_k'$ for $ k \leq n $  in
formula (4.8) can be explicitly computed and are curvatures of the
metric of $\hat g$ of the conformal infinity $(M, \hat g)$. In the
following, we write down the formula for the cases $n =2$ and $
n=4$.

\proclaim{Proposition 4.2} Suppose that $(X^{n+1}, g)$ is a
conformally compact Einstein manifold with even $n$. For a choice
of metric $\hat g\in [\hat g]$ of the conformal infinity $(M,
[\hat g])$, we have
$$
V (X^3, g, \hat g)= - \int_{M^2} \, (\frac {d
}{d s}|_{ s=2} S(s)1 ) \,
dv_{\hat g}, \ \text{for $n=2$} \tag 4.9
$$
and
$$
V (X^5, g, \hat g) = - \frac 1{32\cdot 36}\int_{M^4}
 R_{\hat g} ^2 dv_{\hat g} - \int_{M^4} \,
(\frac {d }{d s}|_{ s=4} S(s)1) \,  dv_{\hat g}
, \ \text{for $n=4$} \tag 4.10
$$
where $R_{\hat g}$ is the scalar curvature for $(M^4, \hat g)$.
\endproclaim

\demo{Proof} For $n=2$, one may calculate and obtain
$$
a_2 (s) = 2 K_{\hat g},
$$
where $K_{\hat g}$ is the Gaussian curvature of $(M^2, \hat g)$, thus $a_2' =0$, and (4.9) follows directly from (4.8).
For $n=4$, one calculate and obtain
$$
a_2(s) =  -\frac {4-s}{4(3-s)}\text{Tr}g^{(2)} ,
$$
thus
$$
a_2' = - \frac 14 \text{Tr} g^{(2)} \tag 4.11
$$
and
$$
\aligned a_4 & = -\frac 1{8(4-s)}((4-s)(2\text{Tr}g^{(4)} -
|g^{(2)}|^2) \\ & \quad \quad - \frac
{(4-s)(6-s)}{4(3-s)}(\text{Tr}g^{(2)})^2 - \frac
{(4-s)}{8(3-s)}\hat\Delta\text{Tr}g^{(2)})\\
& = -\frac 1{8}(2\text{Tr}g^{(4)} - |g^{(2)}|^2 - \frac
{(6-s)}{4(3-s)}(\text{Tr}g^{(2)})^2 - \frac
1{4(3-s)}\hat\Delta\text{Tr}g^{(2)})
\endaligned
$$
where
$$
g^{(2)} = -\frac 12 (Ric_{\hat g} - \frac 16
R_{\hat g} \hat g), \tag 4.12
$$
$$
\text{Tr}g^{(4)} = \frac 14 |g^{(2)}|^2_{\hat g}. \tag 4.13
$$
Thus
$$
a_4' = \frac 1{32}(3(\text{Tr}g^{(2)})^2 + \hat \Delta
\text{Tr}g^{(2)}) \tag 4.14
$$
Insert (4.11) and (4.12) and (4.14) to (4.8) and recall that when
$n=4$, $v^{(2)} = \frac 12 \text{Tr}g^{(2)} $, we obtain (4.10).
\enddemo

We end the section by deriving a variational formula, which indicates that when $n$ is even, the
scattering term in the renormalized volume is a conformal anomaly
and is a conformal primitive of the $Q$-curvature $Q_n$.
Namely,

\proclaim{Theorem 4.3} Suppose that $(X^{n+1}, g)$ is a
conformally compact Einstein manifold with conformal infinity $(M,
[\hat g])$, and that $n$ is even. Then
$$
\frac{d}{d \alpha}|_{\alpha = 0} \int_M \Cal
S_{ e^{ 2 \alpha w} \hat g} \, dv_{e^{2\alpha
w}\hat g} = - 2c_{\frac n2} \int_M w \,
(Q_n)_{\hat g} \, dv_{\hat g} \tag 4.15
$$
where $\Cal S_{\hat g} = \frac d{d s}|_{s=n}
S_{( g, \hat g)}(s)1 $; and $S_{(g, \hat g)}$ is
the scattering operator as defined in (4.3).
\endproclaim

To prove the theorem, we first list some elementary properties of
scattering matrix under conformal change of metrics, we assume the same setting as in Theorem 4.3.

\proclaim{Lemma 4.4} Denote $S(s)=
S^{\hat g} =
S_{(g, \hat g}) (s)$ the scattering
 matrix, and ${\hat g}_w = e^{2w} {\hat g}$ a metric conformal to $\hat g$. Then
\newline \noindent (a)
$$   S_{\hat g_w}(s) = e^{-sw}S(s)e^{(n-s)w}. \tag 4.16 $$
(b) S(s) has a simple pole at $s=n$
and its residue is  $-c_{\frac n2}P_{\frac n2}$, i.e.
$$ S(s) = - \frac {c_{\frac n2}P_{\frac n2}}{s-n} + T(s) \tag 4.17
$$
where where $T(s)$ is the regular part of the scattering matrix
near $s=n$ and $P_{\frac n2}$ is the Paneitz operator.
\newline \noindent (c)
$$ S(n)1 = \lim_{ s \to n } S(s) 1 = T(n) 1 = c_{ \frac n2 } Q_n.
\tag 4.18
$$
\noindent (d)
$$
\aligned e^{n w}S_{\hat g_w}(n) 1 & =
e^{nw}\lim_{s\rightarrow n}
e^{-sw}S(s)e^{(n-s)w} \\
& = \lim_{s\rightarrow n}( - \frac {c_{\frac n2}P_{\frac
n2}e^{(n-s)w}}{s-n} + T(s)e^{(n-s)w}) \\
& = c_{\frac n2}P_{\frac n2}w + S(n)1.
\endaligned
\tag 4.19
$$
\endproclaim

\demo{Proof}\newline
(a) is a simple consequence of the definition of the scattering matrix. (b) and (c) follow from the work of [GZ]. (d) is a consequence of the equation
$$ (P_{\frac n2})_{\hat g}  w + (Q_n)_{\hat g} =
(Q_n)_{\hat g_w} e^{nw} $$
relating the Paneitz operator to $Q$ -curvature on even
dimensional manifolds.
\enddemo

We now compute the variation of
$$
\int_M \frac {d}{ds} |_{ s=n} S_{(g, \hat g)} (s) 1 dv_{\hat g}
$$
under the conformal change of metrics in $[\hat g]$. Denote by $
\Cal S_{\hat g} =  \frac {d}{ds} |_{ s=n}
S_{(g, \hat g)} (s) 1 $ with respect to
the metric $\hat g$ on $M$.

\proclaim{Lemma 4.5}
$$
\int_M  ({\Cal S}_{ e ^{2w}\hat g} e^{nw}  -  \Cal S_{\hat g})
dv_{\hat g} = - c_{\frac n2} \int_M ( w
\,(P_{\frac n2})_{\hat g}  w + 2w (Q_n)_{\hat
g} )  \, dv_{\hat g}. \tag 4.20
$$
\endproclaim

\demo{Proof} By definition we have
$$
{\Cal S}_{ \hat g e^{ 2 w}} e^{ n w} = \lim_{ s \to n} ( \frac {
S_{e^{2w}\hat g} (s) 1 - S_{e^{2w}\hat g} (n) 1 }{ s-n} ) e^{n w}.
\tag 4.21
$$
Apply (4.16) and (4.19), denote $P_{\frac n2} =
(P_{\frac n2})_{\hat g} $, $Q = Q_n =
(Q_n)_{\hat g} $ and $T = T_{(g, \hat g)}$, we obtain
$$
\aligned  e^{n  w}& \frac { S_{e^{2w}\hat g} (s) 1 - S_{e^{2w}\hat
g} (n) 1 }{ s-n} = \frac {e^{(n-s)w}S(s)e^{(n-s) w} - c_{ \frac
n2 } P_{\frac n2}w
- c_{\frac n2 } Q_{\hat g}}{s-n} \\
& = \frac {(e^{(n-s) w}-1)S(s)e^{(n-s) w} + S(s)e^{(n-s) w}- c_{
\frac n2 } P_{\frac n2} w - c_{\frac n2 } Q_{\hat g}}{s-n} \\
& = \frac {e^{(n-s) w}-1}{s-n}S(s)e^{(n-s) w} + \frac{S(s)e^{(n-s)
w}- c_{ \frac n2 } P_{\frac n2}  w - c_{\frac n2 } Q_{\hat
g}}{s-n}.
\endaligned
\tag 4.22
$$

We now claim that after taking limit we have
$$
{\Cal S}_{e^{2 w}\hat g} e^{n w} - {\Cal S}_{\hat g}= - c_{ \frac n2 }
(w  P_{\frac n2}w + wQ_{\hat g}- \frac 12 P_{\frac n2}w^2) - T(n)w. \tag 4.23
$$
To see the claim we observe that
$$
\lim_{ s \to n}\frac {e^{(n-s) w}-1}{s-n}= -w, \tag 4.24
$$
$$
\aligned  \lim_{ s \to n} S(s)e^{(n-s)w} & =  \lim_{ s \to n}(-\frac {c_{ \frac n2 } P_{\frac
n2}e^{(n-s)w}}{s-n} + T(s) e^{(n-s)w}) \\
& = c_{\frac n2}P_{\frac n2}w + c_{\frac n2}Q_n,
\endaligned\tag 4.25
$$
$$
\aligned  & \lim_{s\to n} \frac{S(s)e^{(n-s)w} - c_{\frac
n2}P_{\frac n2}w - c_{\frac n2} Q_{\hat g}}{s-n} \\
= & \lim_{s\to n}( \frac{S(s)(e^{(n-s)w}-1)- c_{ \frac n2 }
P_{\frac n2}w}{s-n} + \frac {S(s)1 - S(n)1}{s-n}) \\
= &  \lim_{s\to n}\frac{S(s)(e^{(n-s)w}-1)- c_{ \frac n2 }
P_{\frac n2}w}{s-n} + \Cal S_{\hat g},
\endaligned \tag 4.26
$$
and
$$
\aligned & \lim_{s\to n} \frac{S(s)(e^{(n-s)w}-1)- c_{ \frac n2 }
P_{\frac n2}w}{s-n} \\
= & \lim_{s\to n}\frac {\frac {-c_{\frac n2}P_{\frac
n2}(e^{(n-s)w}-1)}{s-n} + T(s)(e^{(n-s)w}-1) - c_{\frac
n2}P_{\frac n2}w}{s-n} \\
= & \lim_{s\to n} -c_{\frac n2}P_{\frac n2}\frac {e^{(n-s)w} - 1 -
(n-s)w}{(s-n)^2} + \lim_{s\to n} T(s) \frac {e^{(n-s)w}-1}{s-n} \\
= & -\frac 12 c_{\frac n2}P_{\frac n2}w^2 - T(n) w.
\endaligned
\tag 4.27
$$
Thus the claim (4.23) follows from the formulas (4.24) to (4.27).

Due to the fact that both operators $P_{\frac n2}$ and $T(n)$ are
self-adjoint, we have
$$
\int_M P_{\frac n2}w^2 dv_{\hat g} = \int_M w^2 P_{\frac n2}1
dv_{\hat g} = 0
$$
and
$$
\int_M T(n)w dv_{\hat g} = \int_M wT(n)1 dv_{\hat g} = c_{\frac
n2} \int_M w Q_ndv_{\hat g}.
$$
Thus integrating (4.23), we get
$$
\int_M \Cal S_{e^{2w}\hat g}dv_{e^{2w}\hat g} - \int_M \Cal
S_{\hat g}dv_{\hat g} = -c_{\frac n2}\int_M (wP_{\frac n2}w +
2wQ_n)dv_{\hat g}.
$$
This is the desired formula (4.15).
\enddemo

Theorem 4.3 follows from Lemma 4.5 by a
simple integration.

\vskip .1in
We remark that (4.10) and (4.15) give another proof of the fact
observed in [HS] and [G-1] that when $n=2$, $V(X^3, g)$ is the
conformal primitive of the Gaussian curvature
$K_{\hat g}$ on
$(M^2, \hat g)$; while for $n=4$, $V(X^5,g)$ is the conformal
primitive of $\frac{1}{16} \sigma_2 (\hat g)$ on $(M^4, \hat g)$,
where
$$
\sigma_2 (\hat g) = \frac 16 (R_{\hat
g}^2 - 3|Ric|_{\hat g} ^2)
$$
and the relation of $\sigma_2 (\hat g)$ to $
(Q_4)_{\hat g}$ is given by
$$
(Q_4)_{\hat g} = \frac 16
(- \Delta R + R^2 - 3|Ric|^2)_{\hat g} =
\sigma_2 (\hat g) - \frac 16 \Delta_{\hat g}
R_{\hat g}, \tag 4.28
$$
where $\Delta = \sum \frac {\partial^2}{\partial x_i^2}$ for the
Euclidean metric. We remark that the term $\sigma_2$ plays an
important role in some recent work [CGY] in conformal geometry,
where the sign of $\int_M
(\sigma_2)_{\hat g} dv_{\hat g} $ is used to study existence
of metrics with positive Ricci curvature on compact, closed
manifolds of dimension 4; also the relation (4.28) between
$\sigma_2$ and $Q_4$ plays a crucial role in the proof of the
results in [CGY].

\vskip 0.2in \noindent {\bf References}:

\roster
\item"{[Ab1]}" P. Albin, A renormalized index theorem for some complete
asymptotically regular metrics: the Gauss-Bonnet theorem, preprint
math.DG/0512167.

\vskip 0.1in\item"{[Ab2]}" P. Albin, Renormalizing Curvature
Integrals on Poincare-Einstein Manifolds, preprint math.DG/0504161.

\vskip 0.1in
\item"{[A]}" M. Anderson, $L^2$ curvature and volume renormalization of
the AHE metrics on 4-manifolds, Math. Res. Lett., 8 (2001)
171-188.

\vskip 0.1in
\item"{[Al]}" S. Alexakis, private communication, thesis, Princeton
University, 2003.

\vskip 0.1in
\item"{[B]}" T. Branson, ``Functional determinant", Lecture Note
Series, 4. Seoul National University, 1993.

\vskip 0.1in
\item"{[CGY]}" S.-Y. A. Chang, M. Gursky and P. Yang,
An equation of Monge-Ampere type in conformal geometry and
4-manifolds of positive Ricci curvature, Annals of Math. 155
(2002), 709-787.

\vskip 0.1in
\item"{[CQ]}" S.-Y. A. Chang and J. Qing, The zeta functional
determinants on manifolds with boundary I - The formula, J. Funct.
Analysis, 149 (1997), no.2, 327-362.

\vskip 0.1in
\item"{[CQY]}" S.Y. A. Chang, J. Qing and P. Yang, On the topology
of conformally compact Einstein 4-manifolds,  Noncompact Problems at
the intersection of Geometry, Analysis and Topology; Contempoary Math.
volume 350, 2004, pp 49-61.

\vskip 0.1in
\item"{[E]}" C. Epstein, Appendix in, The divisor of Selberg's Zeta
functions for Kleinian groups, Duke J. Math. 106 (2001) no. 2. 321
- 390.

\vskip 0.1in
\item"{[FG-1]}" C. Fefferman, and R.Graham, Conformal invariants, in
{\it The mathematical heritage of Elie Cartan}, Asterisque, 1985,
95-116.

\vskip 0.1in
\item"{[FG-2]}" C. Fefferman, and R.Graham, Q-curvature and Poincar\'{e}
metrics, Math. Res. Lett.,  9  (2002),  no. 2-3, 139-151.

\vskip .1in
\item"{[FH]}" C. Fefferman and K. Hirachi, Ambient metric constuction of
Q-curvature in conformal and CR geometry, Math. Res.
Lett., 10 (2003), no. 5-6, 819-832.

\vskip 0.1in
\item"{[G-1]}" R. Graham, Volume and Area renormalizations for
conformally compact Einstein metrics,  The Proceedings of the 19th
Winter School "Geometry and Physics" (Srn\`{i}, 1999). Rend. Circ.
Mat. Palermo (2) Suppl. No. 63 (2000), 31--42.

\vskip 0.1in
\item"{[G-2]}" R.Graham, private communication, 2003.

\vskip .1in
\item"{[GL]}"  R. Graham and L. Lee, Einstein metrics with prescribed
conformal infinity on the ball, Adv. Math. 87 (1991), no. 2, 186-225.

\vskip 0.1in
\item"{[GJMS]}" R. Graham, R. Jenne, L. Mason and J. Sparling,
Conformally invariant powers of the Laplacian, I. Existence, J.
London Math. Soc. (2) 46 (1992), 557-565.

\vskip 0.1in
\item "{[GZ]}" C. R. Graham and M. Zworski, Scattering matrix
in conformal geometry, Invent. Math. 152(2003) no. 1, 89-118.

\vskip 0.1in
\item"{[HS]}" M. Henningson and K. Skenderis, The holographic Weyl
anomaly, J. High. Energy Phys. 07 (1998) 023 hep-th/9806087;
Holography and the Weyl anomaly hep-th/9812032.

\endroster

\enddocument